\documentclass{article}

\usepackage{latexsym}                   
\usepackage{epsfig}                     
\usepackage{graphicx}
\usepackage{amssymb}
\usepackage{amsmath}           
\usepackage{amsthm}
\usepackage{enumerate}                  
\usepackage{mathrsfs}                   
\usepackage{algorithm}
\usepackage{algorithmic}
\usepackage{array}
\usepackage{natbib}

\newtheorem{thm}{Theorem}
{\bfseries}{\itshape}
{\bfseries}{\itshape}
{\bfseries}{\itshape}
{\bfseries}{\itshape}
\theoremstyle{definition}
{\bfseries}{\normalfont}
\newtheorem{rem}[thm]{Remark}{\bfseries}{\normalfont}
{\bfseries}{\normalfont}

\newcommand{\bE}{{\mathbb E}}

\newcommand{\bN}{{\mathbb N}}
\newcommand{\bP}{{\mathbb P}}

\newcommand{\cT}{{\mathcal T}}

\newcommand\ring[1]{\mathaccent23{#1}}
\def\Vr{\ring{V}}

\begin{document}
%
%
\thispagestyle{plain}
  \title{Estimating the relative order of speciation or coalescence events on a 
given phylogeny\\
  {\small Running header: Relative order of speciation or coalescence events}}
  
  \author{Tanja Gernhard$^{\rm a,}$\footnote{To whom correspondence
      should be addressed.}, Daniel Ford$^{\rm b,}$\footnote{partially supported by grant NSF-DMS-0241246} , Rutger Vos$^{\rm c}$
and Mike Steel$^{\rm d}$\\
    \date{}
   {\small $^{\rm a}$ Department of Mathematics, Kombinatorische Geometrie (M9), TU M\"{u}nchen } \\ 
   {\small Bolzmannstr. 3, 85747 Garching, Germany} \\ 
   {\small Phone +49 89 289 16882, gernhard@ma.tum.de}\\
 {\small $^{\rm b}$Department of Mathematics, Stanford University, USA}\\
 {\small $^{\rm c}$Department of Biological Sciences, Simon Fraser University, Vancouver, Canada}\\
    {\small $^{\rm d}$Biomathematics Research Centre, University of Canterbury, Christchurch, New Zealand}
  }
  \maketitle

\begin{abstract}
The reconstruction of large phylogenetic trees from data that violates 'clocklike' evolution (or as a supertree constructed from any $m$ input trees)
raises a difficult question for biologists - how can one assign relative dates to the vertices of the tree?  In this paper we investigate this problem, assuming a uniform distribution on the order of the inner vertices of the tree (which includes, but is more general than, the popular Yule distribution on trees). We derive fast algorithms for computing the probability that (i) any given vertex in the tree was the $j$--th speciation event (for each $j$), and (ii) any one given vertex is 'earlier' in the tree than a second given vertex. We show how the first algorithm can be used to calculate the expected length of any given interior edge in any given tree that has been generated under either a constant-rate speciation model, or the coalescent model. 
\end{abstract}


\section{Introduction}

\label{ChaptRank}
A fundamental task in evolutionary biology is constructing evolutionary trees from a variety of data.
These constructed trees show the ancesteral relationship between the species.

Not only the relationship between species is of interest, but also the time between speciation events.
When constructing an evolutionary tree from a set of molecular data which satisfies the molecular clock, the edge lengths can be interpreted as a time scale.
In many cases, no time scale is obtained when constructing a tree though:
\begin{itemize}
\item Often, molecular data does not satisfy the molecular clock and so the edge lengths do not represent a time scale.
\item Trees can be constructed from morphological data or non-standard molecular data like gene order. This does not provide any edge lengths.
\item Having several different trees, one can combine them and construct a `supertree'. Even though there may have been time scales on the original trees, most supertree methods return a tree without a time scale.
\end{itemize}

For those trees, we still want to find edge lengths representing the time between speciation events. In this paper, we will estimate the edge lengths from the shape of the tree. The method works for trees which evolved under the Yule model \citep{Yule1924,Edwards1970,Harding1971,Page1991}. Under the Yule model, in each point of time, each species is equally likely to split. Minor changes to the method for the Yule model give us an edge length estimation for trees under the popular coalescent setting \citep{Nordborg2001}.

\begin{figure}[h]
\begin{center}
\includegraphics[scale = 0.5]{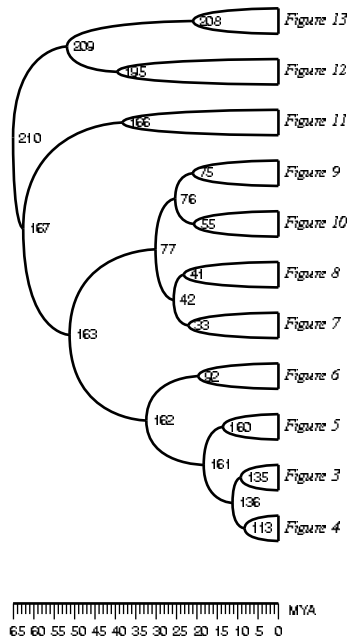}
\caption{Part of the primate supertree. Figure 4 -- 13 are some subtrees, for details see \citep{Vos2006}.}
\label{FigPrimate01}
\end{center}
\end{figure}

An example for a tree with unknown edge lengths is the primate supertree $\cT_p$ recently published in \citep{Vos2006}. Figure \ref{FigPrimate01} shows a part of $\cT_p$. The primate tree is a supertree on 218 species and was constructed with the MRP method (Matrix Representation using Parsimony analysis, see \citep{Baum1992,Ragan1992}).
Since for most of the interior vertices, no molecular estimates were available, the edge lengths for the tree were estimated.
In \citep{Vos2006}, $10^6$ rank functions on $\cT_p$ were drawn uniformly at random. For each of those rank functions, the expected time intervals, i.e. the edge lengths, between vertices were considered (the expected waiting time after the $(n-1)$th event until the $n$th event is $1/n$).
The authors of \citep{Vos2006} concluded their paper by asking for an analytical approach to the estimation of the edge length, which we will provide below.

In order to estimate the edge lengths, we developed the algorithms {\sc RankProb} and {\sc Compare}. Those algorithms answer questions like:

Was speciation event with label 76 in the primate tree (see Fig. \ref{FigPrimate01}) more likely to be an early event in the tree or a late event?
What is the probability that 76 was the 6th speciation event?
Was it more likely that speciation event 76 happened before speciation event 162 or 162 before 76?

The algorithms work for trees where every labeled history is equiprobable. This class of model, which includes the Yule 
model and the coalescent model,
has been popular in macroevolutionary studies \citep{Nee1997,Zhaxybayeva2004}. Note that the algorithms here are the same for the Yule model and the coalescent model, whereas the edge length estimation has minor differences for the two models.

The algorithms {\sc RankProb}, {\sc Compare} and an algorithm for obtaining the expected rank and variance for a vertex were implemented in Python, see \citep{Gernhard2006}.

\section{Probability distribution of the rank of a vertex}
Let $\cT$ be a rooted phylogenetic tree \citep{Steel2003} with $|V|=n$ leaves. The set of interior vertices of $\cT$ shall be $\Vr$. For a binary tree, we have $|\Vr|=n-1$.
Let the function $r$ be a
bijection from the set of interior vertices $\Vr$ of $\cT$ into
$\{1,2, \ldots ,|\Vr|\}$ with $r(v_1) \leq r(v_2)$ if $v_1$ is an ancestor of $v_2$.
The function $r$ is called a {\it rank function} for $\cT$.
A vertex $v$ with $r(v)=i$ is said to have {\it rank} $i$.
Note that $r$ induces a linear order on the set $\Vr$. Further, define $r(\cT):= \{r: \ r \rm{ \ is \ a \ rank \ function \ on} \ \cT \}$.
We are interested in the distribution of the possible ranks for a certain vertex, i.e. we want to know the
probability of $r(v)=i$ for a given~$v \in \Vr$.
If every rank function on a given tree is equally likely, we have
\begin{equation} \label{EqnUnif}
\bP[r(v)=i] = \frac{| \{r: r(v)=i, r \in r(\cT) \}|}{| r(\cT)|}
\end{equation}
which will be calculated for rooted binary trees in polynomial time by algorithm {\sc RankProb}.
In the algorithm, we will use the formula \citep{Steel2003}
\begin{equation} |r(\cT)| = \frac{|\Vr|!}{ \prod_{v \in \Vr} (n_v -1)} \label{EqnNumbRank} \end{equation}
where $n_v$ is the number of leaves below $v$. Note that Equation \ref{EqnNumbRank} holds for binary and nonbinary trees.

\noindent
Examples of stochastic models on phylogenetic trees where each rank function is equally likely include:
\begin{itemize}
\item The Yule model has the probability distribution $\bP[r | \cT] = \frac{\prod_{v \in \Vr} (n_v-1)} {(n-1)!}$ which is the uniform distribution \citep{Edwards1970,Brown1994}.
\item The coalescent model has the same probability distribution on rooted binary ranked trees as the Yule model.
So $\bP[r| \cT]$ is the uniform distribution \citep{Aldous2001}.
\item
For some sets of trees (e.g. those drawn from the uniform model \citep{Pinelis2003}, also known as PDA model), no rank function is induced.  If one assumes that all rank functions are equally likely on these trees, one can apply Equation \ref{EqnUnif} to such trees as well.
\end{itemize}

\subsection{A polynomial-time algorithm}

\begin{figure}
\begin{center}
\input{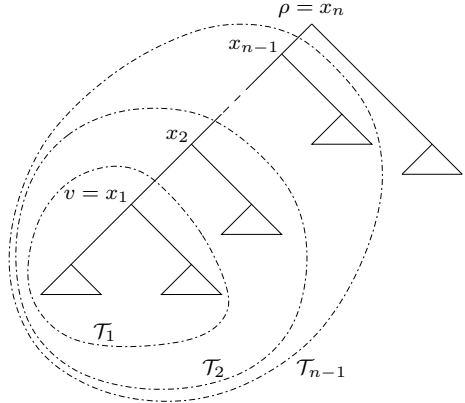}
\caption{Labeling the tree for the algorithm {\sc RankProb}}
\label{FigRankCount1}
\end{center}
\end{figure}

The following algorithm calculates the probability distribution of the rank of a vertex $v$ in a rooted binary phylogenetic tree $\cT$. The idea of the algorithm is the following ({\it cf.} Figure \ref{FigRankCount1}). Label the vertices on the path from $v$ to the root $\rho$ by $v=x_1, \ldots , x_n=\rho$. Let $\cT_m$ be the subtree of $\cT$ containing the vertex $x_m$ and all its descendants. Let $\alpha_{\cT_m,v}(i)$ be the number of rank functions on the tree $\cT_m$ where $v$ has rank $i$. The values ${\alpha}_{\cT_m,v}(i), ~i=1, \ldots, |\Vr|$ are calculated iteratively for $m=1, \ldots, n$. The probability $\bP[r(v)=i]$ equals $\frac{{\alpha}_{\cT_n,v}(i)}{\sum_{i=1}^{|\Vr|}{\alpha}_{\cT_n,v}(i)}$. The $\alpha$-values in the fraction have a lot of factors in common which cancel out. In the following algorithm, we calculate $\alpha$-values without the unnecessary terms instead, $\tilde{\alpha}_{\cT_m,v}(i)$. We have $\alpha_{\cT_m,v}(i) = \tilde{\alpha}_{\cT_m,v}(i) |r(\cT_1)| |r(\cT'_1)| |r(\cT'_2)| \ldots |r(\cT'_{m-1})|$.\\

\noindent
        {\bf Algorithm}: \textsc{RankProb}($\cT,v$)  \index{algorithm R\textsc{ankProb}}\\
        {\bf Input}: A rooted binary phylogenetic tree $\cT$ and an interior vertex $v$.\\
        {\bf Output}: The probabilities $\bP[r(v)=i]$ for ${i=1, \ldots ,|\Vr|}$.
        \begin{algorithmic}[1]
          \STATE Denote the vertices of the path from $v$ to root $\rho$ with \\ $(v=x_1, x_2, \ldots ,x_n=\rho)$.
          \STATE Denote the subtree of $\cT$, consisting of root $x_m$ and all its descendants, by $\cT_m$ for $m=1, \ldots ,n$.
          \STATE Initialize $\tilde{\alpha}_{\cT_m,v}(i):=0$ for $i=1,\ldots,|\Vr_{\cT}|, m=1,\ldots,n$
          \STATE $\tilde{\alpha}_{\cT_1,v}(1):=1$
          \FOR{$m = 2, \ldots ,n$}      
            \STATE $\cT_{m-1}':=\cT_m  \setminus (\cT_{m-1} \cup x_m)$ \qquad ({\it cf.} Figure \ref{FigRankCount2})
            \FOR{$i=m, \ldots ,|\Vr_{\cT_m}|$}
                \STATE $M := \min\{|\Vr_{\cT_{m-1}'}|, i-2\}$
                \STATE $\displaystyle \tilde{\alpha}_{\cT_m,v}(i):=\sum_{j=0}^{M} \tilde{\alpha}_{\cT_{m-1},v} (i-j-1)
                {|\Vr_{\cT_{m-1}}|+|\Vr_{\cT_{m-1}'}|-(i-1) \choose |\Vr_{\cT_{m-1}'}|-j} {i-2 \choose j}  \qquad (\ast)$
            \ENDFOR
          \ENDFOR
          \FOR{$i=1, \ldots , |\Vr_{\cT}|$}
            \STATE $\bP[r(v)=i] := \frac{\tilde{\alpha}_{\cT_n,v}(i)}{\sum_{j} \tilde{\alpha}_{\cT_n,v}(j)}$
          \ENDFOR
          \STATE RETURN $\bP[r(v)=i], i=1,\ldots,|\Vr|.$
        \end{algorithmic}
        \bigskip

\begin{figure}
\begin{center}
\input{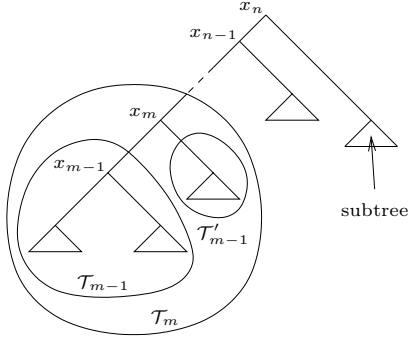}
\caption{Labeling the tree for the recursion in {\sc RankProb}}
\label{FigRankCount2}
\end{center}
\end{figure}

Proving the correctness and runtime of {\sc RankProb} makes use of the following two observations.

\begin{rem} \label{CorSequenceBin}
Let $A_i$ be a set containing $n_i$ elements with a linear order, $i \in \{1,2\}$. There are ${n_1+n_2 \choose n_1}$
possible linear orders on $A_1 \cup A_2$ which preserve the linear order on $A_1$ and $A_2$.
This follows from the observation that the number of such linear orders on $A_1 \cup A_2$ is equivalent to the number of ways of choosing $n_1$ elements from $n_1+n_2$ elements, which is ${n_1+n_2 \choose n_1}$.
\end{rem}

\begin{rem} \label{nchoosek}
The values ${n \choose k}$ for all $n,k \leq N$ ($n,k,N \in \bN$) can be calculated in $O(N^2)$ using Pascal's Triangle.
Thus, after $O(N^2)$ calculations, any value ${n \choose k}$ with $n,k \leq N$ can be obtained
in constant time.
\end{rem}

\begin{thm} \label{ThmRankProb}
{\sc RankProb} returns the quantities $$\bP[r(v)=i]$$
for each given $v\in \Vr$ and all $i \in 1, \ldots ,|\Vr|$. The runtime is $O(|\Vr|^2)$.
\end{thm}

\begin{proof}
Let $\alpha_{\cT_m,v}(i) = \tilde{\alpha}_{\cT_m,v}(i) |r(\cT_1)| |r(\cT'_1)| |r(\cT'_2)| \ldots |r(\cT'_{m-1})|$.
We first show that $\alpha_{\cT_m,v}(i)=|\{r: r(v)=i,~r \in r(\cT_m) \}|$
for $m=1, \ldots ,n,~i=1, \ldots ,|\Vr_{\cT}|$.
That implies $$\bP[r(v)=i] = \frac{| \{r: r(v)=i, r \in r(\cT) \}|}{| r(\cT)|} = \frac{\alpha_{\cT,v}(i)}{\sum_i \alpha_{\cT,v}(i)} = \frac{\tilde{\alpha}_{\cT,v}(i)}{\sum_i \tilde{\alpha}_{\cT,v}(i)}$$ which proves the theorem.

The proof is by induction over $m$.\\
For $m=1$, $\alpha_{\cT_1,v}(1)=|r(\cT_1)| \tilde{\alpha}_{\cT_1,v}(i)= |r(\cT_1)| = |\{r: r(v)=1,~r \in r(\cT) \}|$. Vertex $v$ is the root of $\cT_1$,
so $\alpha_{\cT_1,v}(i)=0$ for all $i>1$.\\
Let $m=k$ and $\alpha_{\cT_m,v}(i)=|\{r: r(v)=i,~r \in r(\cT_m) \}|$
holds for all $m<k$. $\alpha_{\cT_k,v}(i)=0$ clearly holds for all
$i > |\Vr_{\cT_k}|$ since $r_{\cT_k}: v \rightarrow \{1, \ldots
,|\Vr_{\cT_{k}}|\}$. So it remains to verify that the
term $(\ast)$ returns the right values
for $\alpha_{\cT_k,v}(i)$.
Assume that the vertex $v$ is in the $(i-j-1)$-th position in $\cT_{k-1}$ (with $i-j-1>0$) for some rank function $r_{\cT_{k-1}}$ and $v$ shall be in the
$i$-th position in $\cT_k$.

Now combine the linear order in the tree $\cT_{k-1}$ induced by $r_{\cT_{k-1}}$ with a linear order in $\cT_{k-1}'$ induced by $r_{\cT_{k-1}'}$
to get a linear order on $\cT_k$. The first $j$ vertices of $\cT_{k-1}'$ must be inserted between vertices of $\cT_{k-1}$ with lower rank than $v$ so that $v$
ends up to be in the $i$-th position of the tree $\cT_k$.  Count the number of possible way to do this as follows.
The tree $\cT_{k-1}'$ has $|r(\cT_{k-1}')|$ possible rank functions.
Combining a rank function $r_{\cT_{k-1}}$ with a rank function $r_{\cT_{k-1}'}$ to get a rank function $r_{\cT_k}$ with $r_{\cT_k}(v)=i$
means inserting the first $j$ vertices of $\cT_{k-1}'$ anywhere between the first $(i-j-2)$ vertices of $\cT_{k-1}$. There are
$${(i-j-2)+j \choose j} = {i-2 \choose j}$$ possibilities according to Remark \ref{CorSequenceBin}.
For combining the $|\Vr_{\cT_{k-1}}|-(i-j-1)$ vertices of rank bigger than $v$ in $\cT_{k-1}$ with the remaining $|\Vr_{\cT_{k-1}'}|-j$ vertices in $\cT_{k-1}'$,
there are $${|\Vr_{\cT_{k-1}}|-(i-j-1)+|\Vr_{\cT_{k-1}'}|-j \choose |\Vr_{\cT_{k-1}'}|-j} = {|\Vr_{\cT_{k-1}}|+|\Vr_{\cT_{k-1}'}|-(i-1) \choose |\Vr_{\cT_{k-1}'}|-j}$$ possibilities.
This follows again from Remark \ref{CorSequenceBin}.
The number of rank functions $r_{\cT_{k-1}}$ with $r_{\cT_{k-1}}(v)=i-j-1$ is $\alpha_{\cT_{k-1},v} (i-j-1)$ by the induction assumption.
Multiplying all those possibilities gives
$$\alpha_{\cT_{k-1},v} (i-j-1) |r(\cT_{k-1}')| {|\Vr_{\cT_{k-1}}|+|\Vr_{\cT_{k-1}'}|-(i-1) \choose
|\Vr_{\cT_{k-1}'}|-j} {i-2 \choose j}$$
where $\alpha_{\cT_{k-1},v}(i) = \tilde{\alpha}_{\cT_{k-1},v}(i) |r(\cT_1)| |r(\cT'_1)| |r(\cT'_2)| \ldots |r(\cT'_{k-2})|$. The value
$|\{r: r(v)=i,~r \in r(\cT) \}|$ is then the sum over all possible $j$ which establishes the correctness of the algorithm.

All that remains is to verify the runtime. Note that the combinatorial factors ${n \choose k}$ for all $n,k \leq |\Vr|$ can be calculated in advance in quadratic time, see Remark \ref{nchoosek}.
In the algorithm, those factors can then be obtained in constant time.

The most time consuming part of the algorithm is line 13. Adding up all calculations needed for obtaining $\alpha'_{\cT_m,v}(i)$, $m=1, \ldots ,n$, $i=1, \ldots ,|\Vr_{\cT_m}|$ comes to:
$$\sum_{m=2}^n|\Vr_{\cT_m}||\Vr_{\cT_{m-1}'}| \leq \sum_{m=2}^n|\Vr||\Vr_{\cT_{m-1}'}| = |\Vr| \sum_{m=2}^n|\Vr_{\cT_{m-1}'}| \leq |\Vr|^2$$
The last inequality holds since the vertices of the $\cT_{m}'$, $m=1, \ldots ,n-1$, are distinct.
Therefore, the runtime is quadratic.
\end{proof}
\begin{rem} \label{RemExpVar}
With $\bP[r(v)=i]$ from Theorem \ref{ThmRankProb}, the
expected value $\mu_{r(v)}$ and the variance $\sigma_{r(v)}^2$ for
$r(v)$ can be calculated by
\begin{equation}
\mu_{r(v)} = \sum_{i=1}^{|\Vr|} i\bP[r(v)=i] \qquad
\sigma_{r(v)}^2 = \sum_{i=1}^{|\Vr|} i^2 \bP[r(v)=i] - \mu_{r(v)}^2 \notag
\end{equation}
\end{rem}

\begin{rem}
The algorithm {\sc RankProb} can be generalized to non-binary trees \citep{Gernhard2006}. The runtime is again quadratic.
\end{rem}

\section[Estimating edge lengths]{Application of {\sc RankProb} - Estimating edge lengths} \label{EstEdgeLength}
\subsection{The Yule model} \label{Yule}
A very common stochastic model for rooted binary phylogenetic trees with edge lengths is the continuous-time Yule model \citep{Edwards1970}. As in the discrete Yule model, at every point in time, each species is equally likely to split and give birth to two new species. The expected waiting time for the next speciation event in a tree with $n$ leaves is $1/n$. That is, each species at any given time has a constant speciation rate (normalized so that 1 is the expected time until it next speciates).

Assume that the primate tree $\cT_p$ evolved under the continuous-time Yule model.
In \citep{Gernhard2006}, the tree shape of $\cT_p$ (i.e. the tree without edge lengths) under the discrete Yule model is tested against the uniform model and accepts the Yule model.

Here, we describe how to estimate the edge lengths for a tree which is assumed to have evolved under the continuous-time Yule model.

Let $(u,v)$ be an interior edge in $\cT$ with $u$ the immediate ancestor of $v$.
Let $X$ be the random variable `length of the edge $(u,v)$' given that $\cT$ is generated according to the continuous-time Yule model.

\begin{figure}
\begin{center}
\input{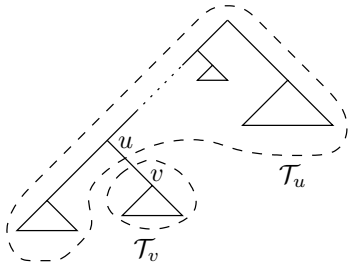}
\caption{Labeling the tree for estimating the edge lengths.}
\label{TreeEstEdge}
\end{center}
\end{figure}

The expected length $\bE[X]$ of the edge $(u,v)$ is given by
$$\bE[X] = \sum_{i,j} \bE[X|r(u)=i, r(v)=j] \bP[r(u)=i, r(v)=j].$$
Since, under the continuous-time Yule model, the expected waiting time for the next speciation event is $1/n$ it follows that:
$$\bE[X|r(u)=i, r(v)=j] = \sum_{k=1}^{j-i} \frac{1}{i+k}.$$
It remains to calculate the probability $\bP[r(u)=i, r(v)=j]$.
This is equivalent to counting all the possible rank functions where $r(u)=i$ and $r(v)=j$.
The subtree $\cT_v$ consists of $v$ and all its descendants.
The tree $\cT_u$ equals the tree $\cT$ where all the descendants of $v$ are deleted, i.e. $v$ is a leaf in $\cT_u$, see Fig. \ref{TreeEstEdge}.

Note that $\bP[r(u)=i, r(v)=j]=0$ if $|\Vr_{\cT_u}|<j-1$.
Therefore, assume $|\Vr_{\cT_u}| \geq j-1$ in the following.

The number of rank functions on $\cT_u$ is $|r(\cT_u)|$.
The probability $\bP[r(u)=i]$ can be calculated with {\sc RankProb}($\cT_u$, $u$).
So the number of rank functions in $\cT_u$ with $\bP[r(u)=i]$ is $\bP[r(u)=i]\cdot |r(\cT_u)|$.

The number of rank functions on $\cT_v$ is $|r(\cT_v)|$.
Let any linear order on the trees $\cT_u$ and $\cT_v$ be given.
Combining those two linear orders into an order, $r$, on $\cT$ with $r(v)=j$ means that the vertices with rank $1,2,\ldots,j-1$ in $\cT_u$ keep their rank. Vertex $v$ gets rank $j$. The remaining $|\Vr_{\cT_u}|-(j-1)$ vertices
in $\cT_u$ and $|\Vr_{\cT_v}|-1$ vertices in $\cT_v$ have to be shuffled together.
According to Remark (\ref{CorSequenceBin}), this can be done in
$$ {|\Vr_{\cT_u}|-(j-1) + |\Vr_{\cT_v}|-1 \choose |\Vr_{\cT_v}|-1} = {|\Vr_{\cT_u}| + |\Vr_{\cT_v}|-j \choose |\Vr_{\cT_v}|-1}$$
different ways. Thus overall there are:
$$\bP[r(u)=i]\cdot |r(\cT_u)| \cdot |r(\cT_v)| \cdot {|\Vr_{\cT_u}| + |\Vr_{\cT_v}|-j \choose |\Vr_{\cT_v}|-1}$$
different rank functions on $\cT$ with $r(u)=i$ and $r(v)=j$.
For the probability $\bP[r(u)=i, r(v)=j]$:
$$\bP[r(u)=i, r(v)=j] = \frac {\bP[r(u)=i]\cdot |r(\cT_u)| \cdot |r(\cT_v)| \cdot {|\Vr_{\cT_u}| + |\Vr_{\cT_v}|-j \choose |\Vr_{\cT_v}|-1}}{\sum_{i,j} \bP[r(u)=i]\cdot |r(\cT_u)| \cdot |r(\cT_v)| \cdot {|\Vr_{\cT_u}| + |\Vr_{\cT_v}|-j \choose |\Vr_{\cT_v}|-1}}$$
Since $|r(\cT_u)|$ and $|r(\cT_v)|$ are independent of $i$ and $j$, those factors cancel out, giving
\begin{equation} \bP[r(u)=i, r(v)=j] = \frac {\bP[r(u)=i]\cdot {|\Vr_{\cT_u}| + |\Vr_{\cT_v}|-j \choose |\Vr_{\cT_v}|-1}}{\sum_{i,j} \bP[r(u)=i] \cdot {|\Vr_{\cT_u}| + |\Vr_{\cT_v}|-j \choose |\Vr_{\cT_v}|-1}}\label{EqnPvu}\end{equation}
Furthermore, note that
$${|\Vr_{\cT_u}| + |\Vr_{\cT_v}|-j \choose |\Vr_{\cT_v}|-1} = \frac{(|\Vr_{\cT}|-j)!}{(|\Vr_{\cT_v}|-1)!(|\Vr_{\cT}|-j-(|\Vr_{\cT_v}|-1))!}$$
Again, since $(|\Vr_{\cT_v}|-1)!$ is independent of $i$ and $j$, this factor cancels out, and so
$$\bP[r(u)=i, r(v)=j] = \frac {\bP[r(u)=i]\cdot \prod_{k=0}^{|\Vr_{\cT_v}|-2} (|\Vr_{\cT}|-j-k) }{\sum_{i,j} \bP[r(u)=i] \cdot \prod_{k=0}^{|\Vr_{\cT_v}|-2} (|\Vr_{\cT}|-j-k)}$$
Let $\Omega = \{(i,j): i<j, i,j \in \{1, \ldots, |\Vr| \}, |\Vr_{\cT_u}| \geq j-1\}$. With this notation, the expected edge length $\bE[X]$ is
\begin{eqnarray}
\bE[X] &=& \sum_{(i,j) \in \Omega}  \bE[X|r(u)=i, r(v)=j] \bP[r(u)=i, r(v)=j] \notag \\
&=& \sum_{(i,j) \in \Omega} \left[ \left( \sum_{k=1}^{j-i} \frac{1}{i+k} \right) \frac {\bP[r(u)=i]\cdot \prod_{k=0}^{|\Vr_{\cT_v}|-2} (|\Vr_{\cT}|-j-k) }{\sum_{(i,j) \in \Omega} \left[ \bP[r(u)=i] \cdot \prod_{k=0}^{|\Vr_{\cT_v}|-2} (|\Vr_{\cT}|-j-k) \right]} \right] \notag \\
&=& \frac {\sum_{(i,j) \in \Omega} \left[ \left( \sum_{k=1}^{j-i} \frac{1}{i+k} \right) \cdot \bP[r(u)=i]\cdot \prod_{k=0}^{|\Vr_{\cT_v}|-2} (|\Vr_{\cT}|-j-k)  \right]}{\sum_{(i,j) \in \Omega} \left[ \bP[r(u)=i] \cdot \prod_{k=0}^{|\Vr_{\cT_v}|-2} (|\Vr_{\cT}|-j-k) \right]}  \label{EqnExpEdgeLength}
\end{eqnarray}

\bigskip

\begin{rem}
Equation \ref{EqnExpEdgeLength} enables the estimation of the length of every interior edge. For pendant edges, the approach above gives no definite answer. All we know is that the time from the latest interior vertex, which has rank $n-1$, until today is expected to be at most $1/n$ where $n$ is the number of leaves.

Suppose that the growth process is stopped as soon as the $n-1$-st speciation event occurs. In this case the expected length $X$ of a pendant edge below an interior vertex $v$ is: \begin{eqnarray} \bE[X] &=& \sum_{i=1}^{n-1} \bP[r(v)=i] \sum_{k=i}^{n-2} \frac{1}{k+1} \notag \end{eqnarray} The expected depth of vertex $v$ from the first branchpoint is: \begin{equation} \sum_{i=1}^{n-1} \bP[r(v)=i] \sum_{k=1}^{i-1} \frac{1}{k+1} \notag \end{equation} So the depth $Y$ of the leaf in question from the first branchpoint has expectation independent of $v$: \begin{eqnarray} \bE[Y] &=& \sum_{i=1}^{n-1} \bP[r(v)=i] \sum_{k=1}^{i-1} \frac{1}{k+1} + \sum_{i=1}^{n-1} \bP[r(v)=i] \sum_{k=i}^{n-2} \frac{1}{k+1} \notag \\ &=& \sum_{i=1}^{n-1} \bP[r(v)=i] \sum_{k=1}^{n-2} \frac{1}{k+1} \notag\\ &=& \sum_{k=1}^{n-2} \frac{1}{k+1} \notag \end{eqnarray} In other words, assigning to each edge of a given tree topology its expected length gives a tree which obeys a molecular clock.

\end{rem}

\begin{rem}
Often, an inferred tree has vertices with more than two descendants, i.e. there is lack of resolution due to, e.g. confliciting data.
Our calculation for the expected edge length assumes a binary tree though.

However, the expected edge length may be calculated for each possible binary resolution of the supertree.
Assume the supertree $\cT$ has the possible binary resolutions $\cT_1, \ldots, \cT_m$.
For an edge $(u,v)$ in $\cT$ where $u$ is the immediate ancestor of $v$, the expected edge length is calculated in the trees $\cT_i$ for $i = 1, \ldots, m$. The expected edge length in $\cT_i$ is denoted by $e_i$ for $i = 1, \ldots, m$.
Note that if $u$ is a vertex with more than two descendants in $\cT$ then $v$ is in general not a direct descendant of $u$ in $\cT_i$. The value $e_i$ in resolution $\cT_i$ is then the sum of all expected edge lengths on the path from $u$ to $v$ in $\cT_i$.

Calculate the expected edge length $\bE[X]$ of $(u,v)$ in the supertree $\cT$ by
\begin{equation}
\bE[X]=\frac{\sum_i e_i \bP[\cT_i]}{\sum_i \bP[\cT_i]} \label{EqnWeightSum}
\end{equation}
where the probability of a tree $\cT$ under the Yule model is \citep{Brown1994} $$\bP[\cT] = \frac{2^{n-1}} {n! \prod_{v \in \Vr} (n_v-1)}$$

Again, once the expected length of pendant edges is included the resulting tree obeys a molecular clock, meaning that all leaves are at the same depth.
\end{rem}

\subsection{The coalescent process} \label{coalescent}
The edge length estimation in the previous section works for the continuous-time Yule model. By changing the method above slightly, we get an edge length estimation for the coalescent process. In the coalescent setting, we have
$$\bE[X|r(u)=i, r(v)=j] = \sum_{k=1}^{j-i} \frac{1}{(i+k)(i+k-1)}.$$
Therefore, the expected edge length for an interior edge $(u,v)$ can be calculated by the following modification of Equation \ref{EqnExpEdgeLength}:

\begin{eqnarray}
\bE[X]
= \frac {\sum_{(i,j) \in \Omega} \left[ \left( \sum_{k=1}^{j-i} \frac{1}{(i+k)(i+k-1)} \right) \cdot \bP[r(u)=i]\cdot \prod_{k=0}^{|\Vr_{\cT_v}|-2} (|\Vr_{\cT}|-j-k)  \right]}{\sum_{(i,j) \in \Omega} \left[ \bP[r(u)=i] \cdot \prod_{k=0}^{|\Vr_{\cT_v}|-2} (|\Vr_{\cT}|-j-k) \right]} \notag
\end{eqnarray}

\bigskip

The calculations in Section \ref{Yule} and \ref{coalescent} provide exact values for the expected length of an interior edge under the Yule or coalescent process as an alternative to simulations. However simulations also provide some indication of the variability in the estimate of edge lengths, and it may be of interest to also investigate analytically the variance (or even the distribution) of the edge length in future work, rather than just its mean. 

\section{Comparing two interior vertices}
The algorithm {\sc RankProb} can also be used for comparing two interior vertices. Assume again that every rank function on a rooted binary phylogenetic tree $\cT$ is equally likely.
The aim is to compare two interior vertices $u$ and $v$ of $\cT$. Was $u$ more likely before (of lower rank than) $v$ or $v$ before $u$? 
In other words, what is the probability
\label{ProbComp}
$$\bP_{u<v}:=\bP[r(u) < r(v)]$$
where $r(T)$ is the set of all possible rank functions on $\cT$.
Note that it does not hold
$\bP[r(u)<r(v)] = \bP[r(u)>r(v)]$ even with the uniform distibution on the rank functions.
The probability $\bP_{u<v}$ is equivalent to counting all the possible rank functions on $\cT$ in
which $u$ has lower rank than $v$ and divide that number by all possible rank functions on $\cT$.
One idea is to sum up the probabilities $\bP[r(u)=i, r(v)=j]$ in Equation \ref{EqnPvu} for all $i<j$ which yields to a runtime of $O(|V|^4)$.
The following algorithm {\sc Compare} solves the problem in quadratic time. In the following, for a vertex $v$, the subtree $\cT_v$ of $\cT$ consists again of $v$ and all its descendants.\\

\noindent
        {\bf Algorithm} {\sc Compare} ($\cT,u,v$) \index{algorithm C{\sc ompare}}\\
        {\bf Input}: A rooted binary phylogenetic tree $\cT$ and two distinct interior vertices $u$ and $v$.\\
        {\bf Output}: The probability $\bP_{u<v} :=\bP[r(u) < r(v) | \cT]$.
        \begin{algorithmic}[1]
          \STATE Denote the most recent common ancestor of $u$ and $v$ by $\rho_1$.
          \IF {$\rho_1 = v$}
            \STATE RETURN $\bP_{u<v} = 0$.
          \ENDIF
          \IF {$\rho_1 = u$}
            \STATE RETURN $\bP_{u<v} = 1$.
          \ENDIF
          \STATE Let $\cT_{\rho_1}$ be the subtree of $\cT$ which is induced by $\rho_1$.
          \STATE Delete the vertex $\rho_1$ from $\cT_{\rho_1}$. The two evolving subtrees are labeled $\cT_u$ and $\cT_v$ with $u \in \cT_u$ and $v \in \cT_v$.
          \STATE Run {\sc RankProb($\cT_u,u$)} and {\sc RankProb($\cT_v,v$)} to get $\bP[r(u)=i]$ on $\cT_u$ and $\bP[r(v)=i]$ on $\cT_v$ for all possible $i$.
            \FOR{$i = 1, \ldots ,|\Vr_{\cT_u}|$}
                \STATE $ucum(i) := \sum_{k=1}^i \bP[r(u)=i]$
            \ENDFOR
                \STATE $\bP_{u<v} := 0$
                \FOR{$i = 1, \ldots, |\Vr_{\cT_v}|$}
                \FOR{$j = 1, \ldots |\Vr_{\cT_u}|$}
                \STATE $p := \bP[r(v)=i] \cdot {i-1+j \choose j} \cdot {|\Vr_{\cT_v}|-i+|\Vr_{\cT_u}|-j \choose |\Vr_{\cT_u}|-j} \cdot ucum(j)$
                    \STATE $\bP_{u<v} := \bP_{u<v}+p$
                \ENDFOR
            \ENDFOR
                \STATE $tot := {|\Vr_{\cT_u}|+|\Vr_{\cT_v}| \choose |\Vr_{\cT_v}|}$
                \STATE $\bP_{u<v} := \bP_{u<v}/tot$
                \STATE RETURN $\bP_{u<v}$
          \end{algorithmic}
        \bigskip

\begin{thm}
The algorithm {\sc Compare} returns the value $$\bP_{u<v} =\bP[r(u) < r(v)].$$ The runtime of {\sc Compare} is $O(|\Vr|^2)$.
\end{thm}
\begin{proof}
Note that the probability of $u$ having smaller rank than $v$ in tree $\cT_{\rho_1}$ equals the probability of $u$ having smaller
rank than $v$ in tree $\cT$, since for any rank function on $\cT_{\rho_1}$, there is the same number of linear extensions to get a rank function on the tree $\cT$.

So it is sufficient to calculate the probability $\bP_{u<v}$ in $\cT_{\rho_1}$.
If $\rho_1=u$ then $u$ is an ancestor of $v$ in $\cT$, so return $\bP_{u<v}=1$.
If $\rho_1=v$ then $v$ is an ancestor of $u$ in $\cT$, so return $\bP_{u<v}=0$.

Now assume that $\rho_1 \neq u$ and $\rho_1 \neq v$.
The run of {\sc RankProb} calculates the probability $\bP[r(u)=i]$ in the tree $\cT_u$ and $\bP[r(v)=i]$ in $\cT_v$ for all $i$.
Next, combine those two linear orders. Assume that $r(v)=i$ and that $j$ vertices of $\cT_u$ are inserted before $v$.
Inserting $j$ vertices of $\cT_u$ into the linear order of $\cT_v$ before $v$ is possible in
${i-1+j \choose j}$ ways (see Remark \ref{CorSequenceBin}).
Putting the remaining vertices in a linear order is possible in
${|\Vr_{\cT_v}|-i+|\Vr_{\cT_u}|-j \choose |\Vr_{\cT_u}|-j}$ ways.
The probability that the vertex $u$ is among the $j$ vertices which have smaller rank than $v$ is $\bP[r(u) \leq j] = ucum(j)$. There are $|r(\cT_u)|$ possible linear orders on $\cT_u$ and $|r(\cT_v)|$ possible linear orders on $\cT_v$.
The number of linear orders where vertex $v$ has rank $i$ in $\cT_v$, $v$ has rank $i+j$ in $\cT_{\rho_1}$ and $r(u)<i+j$ therefore equals
$$p'_{i,j} = \bP[r(v)=i] \cdot |r(\cT_v)| \cdot {i-1+j \choose j} \cdot {|\Vr_{\cT_v}|-i+|\Vr_{\cT_u}|-j \choose |\Vr_{\cT_u}|-j} \cdot ucum(j) \cdot |r(\cT_u)|$$
Adding up the $p'$ for each $i$ and $j$ gives the number of linear orders where $u$ has smaller rank than $v$.

Combining a linear order on $\cT_v$ with a linear order on $\cT_u$ is possible in
$$tot := {|\Vr_{\cT_u}|+|\Vr_{\cT_v}| \choose |\Vr_{\cT_v}|}$$
different ways (see Remark \ref{CorSequenceBin}). There are $|r(\cT_u)|$ linear orders on $\cT_u$ and $|r(\cT_v)|$ linear orders on $\cT_v$, so on $\cT_{\rho_1}$, there are
$$tot' := {|\Vr_{\cT_u}|+|\Vr_{\cT_v}| \choose |\Vr_{\cT_v}|} |r(\cT_v)| |r(\cT_v)|$$
linear orders.
Therefore:
$$\bP_{u<v} = \frac{\sum_{i,j} p'_{i,j}}{tot'} = \frac{\sum_{i,j} p_{i,j}}{tot}$$
with $p_{i,j} = \bP[r(v)=i] \cdot {i-1+j \choose j} \cdot {|\Vr_{\cT_v}|-i+|\Vr_{\cT_u}|-j \choose |\Vr_{\cT_u}|-j} \cdot ucum(j)$. This shows that {\sc Compare} works correct.

Since {\sc RankProb} has quadratic runtime, {\sc Compare} also has quadratic runtime.
\end{proof}

\section{Acknowledgements}
We thank Arne Mooers for very helpful comments and suggestions on earlier versions of this manuscript  and the two anonymous referees for a very careful report.

\nocite{Hey1992}
\bibliographystyle{abbrvnat}
\bibliography{bibliography1}
\end{document}